\documentclass{birkjour}
\usepackage{subfigure}
\usepackage{color}
\usepackage{verbatim}
 \newtheorem{thm}{Theorem}[section]
 \newtheorem{corollary}[thm]{Corollary}
 \newtheorem{lemma}[thm]{Lemma}
 \newtheorem{Proposition}[thm]{Proposition}
 \theoremstyle{definition}
 
 \theoremstyle{remark}
 
 \newtheorem{example}{Example}
 \numberwithin{equation}{section}

 \newcommand{\R}{\mathbb{R}}

\renewcommand{\L}{\mathcal{L}}

\begin{document}

%
%

\title[Line Congruences]
 {{\LARGE Singularities of Generic Line Congruences}}

\author[M.Craizer]{Marcos Craizer}

\address{%
Departamento de Matem\'{a}tica- PUC-Rio\br
Rio de Janeiro, RJ, Brasil}
\email{craizer@puc-rio.br}

\author[R.A.Garcia]{Ronaldo Garcia}

\address{%
Instituto de Matem\'atica e Estat\'istica- UFG\br
Goi\^ania, GO, Brasil}
\email{ragarcia@ufg.br}

\thanks{The authors want to thank CNPq and CAPES (Finance Code 001)  for financial support during the preparation of this manuscript. \newline E-mail of the corresponding author: craizer@puc-rio.br}



\date{April 13, 2020}

\begin{abstract}
Line congruences are $2$-dimensional families of lines in $3$-space.  The singularities that appear in generic line congruences
are folds, cusps and swallowtails (\cite{Izu}). In this paper we give a geometric description of these singularities. The main tool used is the existence of an equiaffine pair defining a generic line congruence.
\end{abstract}

\maketitle

\section{Introduction}

Line congruences are $2$-dimensional families of lines in $3$-space. This is a classical topic of projective differential geometry (\cite{Lane}). The most popular example of a line congruence is the family of euclidean normal lines of a surface, but we can consider also the Blaschke  normals of a non-degenerate surface, normals of a convex surface with respect to a non-euclidean norm in $\mathbb{R}^3$, lines connecting points of a pair of surface patches with parallel tangents, among other examples. 

In this paper we shall be interested only in local questions concerning the line congruence, and so we consider the parameter space as a subset $U\subset\mathbb{R}^2$. For each $(u,v)\in U$, we denote by $f(u,v)$ a base point and by $\xi(u,v)$ a vector
in the direction of the corresponding line. Thus the line congruence can be written as a map $F:U\times\mathbb{R}\to\mathbb{R}^3$ given by 
\begin{equation}\label{eq:CongruenceMap}
F(u,v,t)=f(u,v)+t\xi(u,v).
\end{equation}
In \cite{Izu}, it is proved that the generic singularities of such map are folds, cuspidal edges and swallowtails. In this paper, we shall describe geometrically the precise conditions for their occurrences.

Line congruences determine {\it principal directions} as solutions of a binary differential equation (BDE). Denoting by $\delta$ the {\it discriminant} of the BDE, the set $\delta=0$ is generically a smooth curve dividing the parameter space in regions with two principal directions and no principal directions. At the {\it discriminant curve} $\delta=0$ there is only one principal direction. Points of the discriminant curve with tranversal principal direction are called {\it regular}, while points with tangent principal direction are called {\it singular}.

Outside the discriminant curve, the singularities of a line congruence occur at certain points called {\it ridge points}. 
Generically, these points also form a regular curve in the parameter space. We prove that at ridge points with transversal principal direction, the singularity is a cuspidal edge, while at points with tangent principal directions, the singularity is generically a swallowtail. Altgough these results can obtained from a more general result (Prop.1.3. of \cite{Kokubu}), our approach using jets has the advantage of providing explicit examples of each type of singular point. 

Returning to the discriminant curve, we give a geometric characterization of the generic singularities. We show 
that, at a regular point, the singularity is a fold, while at a singular point, the singularity is generically a cuspidal edge. Moreover,
the ridge curve touches the discriminant curve only at singular points and the condition for a singular point 
to be a cuspidal edge is that the contact between the ridge curve and the separatrix is of order two.

In the proof of the results, we shall consider the {\it support function} of an {\it equiaffine pair}. It happens that the {\it bifurcation set} of the support function of an equiaffine pair $(f,\xi)$ coincides with the {\it focal set} of the line congruence,
which is the image of the singular set of the map \eqref{eq:CongruenceMap}. 
Our main tool is the existence, outside umbilical points, of an equiaffine pair representing an arbitrary line congruence.
Since umbilical points do not occur in generic line congruences, we can use equiafine pairs to prove our results.

The paper is organized as follows: In section 2 we define the basic elements of a line congruence, in section 3 we study the ridge curves outside the discriminant set and in section 4 we descibe the geometry of the principal directions and ridge curves at the discriminant set. In section 5 we discuss properties of equiaffine pairs and in section 6 we prove that, outside umbilical points, every line congruence can be represented by an equiaffine pair. In section 7, respectively 8, we use the equiaffine pairs to describe geometrically the singularities of the line congruence ouside the discriminant set,  respectively,  at the discriminant set.

\newpage

\section{Elements of a Line Congruence}\label{sec:LocalModel}


\subsection{Space of Line Congruences}

In this section we use homogeneous coordinates.  

\paragraph{Pl\"ucker coordinates}
Consider a line defined by $2$ points $y=\left[ y_1:y_2:y_3:y_4 \right]$ and  $z=\left[ z_1:z_2:z_3:z_4 \right]$ in the projective $3$-space $\mathbb{P}^3$ and let 
\begin{equation}\label{eq:DefineOmega}
\omega_{ij}=y_iz_j-y_jz_i, 
\end{equation}
$1\leq i,j \leq 4$, $i\neq j$. The Pl\"ucker (homogeneous) coordinates of the line is defined by 
\begin{equation*}
\omega=\left[  \omega_{12}: \omega_{34}:  \omega_{13}: \omega_{42}: \omega_{14}: \omega_{23}\right]
\end{equation*}
in the projective space $\mathbb{P}^5$. Then $\omega$ belongs to the hyperquadric $Q$ defined by
\begin{equation}\label{eq:Quadric}
\omega_{12}\omega_{34}+\omega_{13}\omega_{42}+\omega_{14}\omega_{23}=0.
\end{equation}
The space of lines is thus represented by a $4$-dimensional submanifold $Q$ of $\mathbb{P}^5$ (for details see \cite{Lane}, p. 72). 

Consider a line $q_0$ in the affine plane $\omega_{34}=-1$. In a neighborhood of $q_0$, Equation \eqref{eq:Quadric} van be written as
\begin{equation}\label{eq:Quadric1}
\omega_{12}=\omega_{13}\omega_{42}+\omega_{14}\omega_{23}.
\end{equation}
Thus, in this neighborhood, a line is represented by the four independent coordinates 
$$
\bar\omega=\left( \omega_{13}, \omega_{42}, \omega_{14}, \omega_{23}\right).
$$

\smallskip

\paragraph{Immersive line congruences} Consider a line congruence $L$ defined in the domain $U$ and assume that $(0,0)\in U$ with $L(0,0)=q_0$ in the affine plane
$\omega_{34}=-1$. Then, by the above considerations, we can consider the line congruence as a smooth map 
$L:U\subset\mathbb{R}^2\to\mathbb{R}^4$. We shall denote by 
$\mathcal{U}$ the space line congruences and by $\mathcal{U}_{im}\subset\mathcal{U}$ the subspace of maps 
$L:U\to\mathbb{R}^4$ that are immersions. It is not difficult to prove that $\mathcal{U}_{im}$ 
is open and dense in $\mathcal{U}$ (see for example \cite{Bruce-Giblin}, prop.8.19). 

\smallskip

\paragraph{Transversal surfaces}
Denote by $\mathcal{U}_0$ the space of line congruences $L$ defined in a neighborhood $U$ of $(0,0)$ such that there exists a smooth regular surface $f(u,v)$ that is transversal to the lines $L(u,v)$. It is easy to see that a line congruence $L$ belongs to $\mathcal{U}_0$
if and only if it can be written in the form \eqref{eq:CongruenceMap} with 
\begin{equation}\label{model}
f(u,v)=(u,v,0),\ \ \ \xi(u,v)=(\bar\xi(u,v),1), \ \ (u,v)\in U,
\end{equation}
where $U$ is some neighborhood of $(0,0)$ and $\bar\xi=(\xi_1,\xi_2):U\to\R^2$ is some smooth map.

\begin{Proposition}
The space $\mathcal{U}_0$ is open and dense in $\mathcal{U}_{im}$.
\end{Proposition}

\begin{proof}

Consider a line congruence $L\in\mathcal{U}_{im}$. Then one of the following determinants must be non-zero
\begin{equation*}
\left|
\begin{array}{cc}
\frac{\omega_{14}}{\partial u} & \frac{\omega_{42}}{\partial u} \\
\frac{\omega_{14}}{\partial v} & \frac{\omega_{42}}{\partial v}
\end{array}
\right|
\ \ , \ \
\left|
\begin{array}{cc}
\frac{\omega_{13}}{\partial u} & \frac{\omega_{23}}{\partial u} \\
\frac{\omega_{13}}{\partial v} & \frac{\omega_{23}}{\partial v}
\end{array}
\right|
\ \ , \ \ 
\left|
\begin{array}{cc}
\frac{\omega_{23}}{\partial u} & \frac{\omega_{42}}{\partial u} \\
\frac{\omega_{23}}{\partial v} & \frac{\omega_{42}}{\partial v}
\end{array}
\right|
\ .
\end{equation*}

\smallskip\noindent
{\bf Case 1:} In the first case, the line congruence can be seen 
as a map 
$
(\omega_{14},\omega_{42})=L(\omega_{13},\omega_{23}).
$
Since $\omega_{34}=-1$, the lines are generated by two points of the form 
$$
\left[ u:v:0:1\right],\ \ \left[ p:q:1:1\right],
$$
which implies that
$$
\omega_{14}=u-p,\ \omega_{42}=q-v,\ \omega_{13}=u,\ \omega_{23}=v.
$$
We conclude that $(p,q)$ are functions of $(u,v)$, which proves the proposition in this case.

\smallskip\noindent
{\bf Case 2:} The second case is similar. The line congruence can be seen 
as a map 
$(\omega_{13},\omega_{23})=L(\omega_{14},\omega_{42}),$
and the lines are defined by points of the form
$$
\left[ p:q:1:1\right],\ \ \left[ u:v:1:0\right],
$$
which implies that
$$
\omega_{14}=-u,\ \omega_{42}=v,\ \omega_{13}=p-u,\ \omega_{23}=q-v.
$$
We conclude that $(p,q)$ are functions of $(u,v)$, which proves the proposition in this case.

\smallskip\noindent
{\bf Case 3:} In the third case, the line congruence can be seen 
as a map 
$(\omega_{23},\omega_{42})=L(\omega_{13},\omega_{14}).$
Choose points in the lines of the form
$$
\left[ u:p:0:1\right],\ \ \left[ v:q:1:1\right],
$$
which implies that
$$
\omega_{14}=u-v,\ \omega_{42}=q-p,\ \omega_{13}=u,\ \omega_{23}=p,
$$
and we conclude that $p$ and $q$ are functions of $(u,v)$. If $p_v(0,0)\neq 0$, we can make the change the variables $U=u,V=p(u,v)$ to obtain our result. Similarly if $q_u(0,0)\neq 0$ we can make the change the variables $U=q(u,v), v=V$ to obtain the result. 
Assume then that $q_u=p_v=0$ at the origin. 

By Thom's transversality theorem, the complement of the set of line congruences satisfying
the three equations $q_u=p_v=0$ and $q_v=p_u$ is open and dense. Thus we can consider only line congruences as above satisfying the property that $q_v\neq p_u$ at the origin.
For such line congruences, the map
$$
(u,v)\to \left[ \frac{1}{2}(u+v): \frac{1}{2}(p+q):\frac{1}{2}:1 \right]
$$
is an immersion transversal to the lines $L(u,v)$, which proves the proposition.

\end{proof}

The space $\mathcal{U}_0$ can be identified with the space of smooth maps $\bar\xi:U\to\R^2$ endowed with the $C^{\infty}$-topology. From now on we shall consider only line congruences in $\mathcal{U}_0$.

\subsection{Principal directions and focal set}

Any $1$-parameter family of lines $(u(t),v(t))$ determines a ruled surface. This ruled surface is developable
if and only if 
$\left[  \xi, f_t,\xi_t \right]=0,$
which is equivalent to 
\begin{equation*}\label{eq:CurvatureCongruence1}
\left[ \xi, f_u,\xi_u\right]du^2+(\left[ \xi, f_v,\xi_u\right]+\left[ \xi, f_u,\xi_v\right])dudv+\left[ \xi, f_v,\xi_v\right]dv^2=0.
\end{equation*}
The solutions of the above Binary Differential Equation (BDE) are called the {\it principal (curvature) lines} of the congruence.
Writing
\begin{equation}\label{eq:DefineShape}
\left\{
\begin{array}{c}
\xi_u=-af_u-cf_v+\tau_1\xi\\
\xi_v=-bf_u-df_v+\tau_2\xi
\end{array}
\right.
\end{equation}
the above BDE becomes
\begin{equation}\label{eq:CurvatureCongruence2}
cdu^2+(d-a)dudv-bdv^2=0.
\end{equation}
Equation \eqref{eq:DefineShape} can be written more compactly in the form 
\begin{equation}\label{eq:DefineShape2}
D_X\xi=-f_*(SX)+\tau(X)\xi,
\end{equation}
where $X\in\mathfrak{X}(U)$, $\tau=(\tau_1,\tau_2)$ is a $1$-form in $U$ and
$$
S=
\left[
\begin{array}{cc}
a & b\\
c & d
\end{array}
\right]
$$
is a linear map called {\it shape operator}. Note that the eigenvectors of the shape operator 
are tangent to the principal lines of the congruence.
If we change $f$ by $f+\lambda\xi$, then the new shape operator is $S(I-\lambda S)^{-1}$, and if we change $\xi$ by $\lambda\xi$, $S$ is multiplied by $\lambda$. In any case, the eigenvectors of $S$ are not changed. Thus we can say that the principal lines of the congruence is independent of the choice of $f$ or $\xi$.


The {\it focal set} is the image of the singular set of the map \eqref{eq:CongruenceMap}.
Since, for $X\in\mathfrak{X}(U)$, $t\in\mathbb{R}$,
$$
D_XF= f_*((I-tS)X)+t\tau(X)\xi, \ \ D_tF=\xi,
$$
we conclude that $F$ is singular if and only if $t=\lambda_i^{-1}$, where $\lambda_i$,  $i=1,2$, denote the eigenvalues of the shape operator. Thus the focal set is the union $F_1\cup F_2$, where
\begin{equation}\label{eq:FocalSurface}
F_i=f+\lambda_i^{-1}\xi.
\end{equation}
The focal set is the envelope of the lines of the line congruence and is thus independent of the choice of $f$ and $\xi$.

\subsection{Ridge and discriminant sets}

We say that $(u,v)$ is a $i$-ridge point if the focal surface $F_i$ is singular at this point. 
Differentiating Equation \eqref{eq:FocalSurface} we obtain,  for $X\in\mathfrak{X}(U)$,
\begin{equation}\label{eq:Ft}
D_XF_i=f_*\left((I-\lambda_i^{-1}S)X\right)+\left( \tau(X)\lambda_i^{-1}-d\lambda_i(X)\lambda_i^{-2}\right)\xi.
\end{equation}
Thus
$D_XF_i=0$ if and only if $X$ is a principal direction associated with the eigenvalue $\lambda_i$ and 
$$
d\lambda_i(X)-\lambda_i\tau(X)=0.
$$
One can verify without difficulties that the ridge set is independent of the choice of $f$ and $\xi$. 
We shall verify below that generically the ridge set is a smooth curve.


The eigenvalues $\lambda_i$, $i=1,2$, of $S$ are the solutions of the quadratic equation 
$\lambda^2-(a+d)\lambda+ad-bc=0.$ Denote by
\begin{equation}\label{eq:Discriminant}
\delta(u,v)=(a-d)^2+4bc
\end{equation}
the discriminant of this equation. Then there are $2$ distinct eigenvectors in the set $\delta>0$ and no eigenvectors in the set $\delta<0$. In the set $\delta=0$, the {\it discriminant set}, the principal directions are coincident. It is clear that the set $\delta=0$ depends only on the congruence and not on the choice of $f$ and $\xi$. 

In order to guarantee that $\delta=0$ is a regular curve, we must verify that $0$ is a regular value of $\delta$. 
Denote by $\mathcal{U}_1\subset\mathcal{U}_0$ the set of line congruences where this occurs. 
A line congruence belongs to $\mathcal{U}_0\setminus\mathcal{U}_1$ if and only if $\delta=\delta_u=\delta_v=0$. 
It is clear that line congruences in $\mathcal{U}_1$ do not admit umbilical points.

The condition for the regularity of $\delta$ involves $3$ equations in $2$ variables, and the equations involve up to second derivatives of $\bar{\xi}$. Thus, by
Thom's transversality theorem, $\mathcal{U}_1\subset\mathcal{U}_0$ is open and dense. From now on we shall always assume that the line congruence belongs to $\mathcal{U}_1$.

\subsection{Local model and notation}

Let 
$$
f(u,v)=\left(u,v,0 \right),\ \ \xi(u,v)=\left(\xi_1(u,v),\xi_2(u,v),1\right)
$$ 
and write
\begin{equation}
\left\{
\begin{array}{c}
\xi_u=-af_u-cf_v\\
\xi_v=-bf_u-df_v
\end{array}
\right.
\end{equation}
The $k$-jets of $a,b,c$ and $d$ are denoted
$$
a=a_{0}+a_{10}u+a_{01}v+\frac{1}{2}\left( a_{20}u^2+2a_{11}uv+a_{02}v^2 \right)+\cdots,
$$
$$
b=b_{0}+b_{10}u+b_{01}v+\frac{1}{2}\left( b_{20}u^2+2b_{11}uv+b_{02}v^2 \right)+\cdots,
$$
$$
c=c_{0}+c_{10}u+c_{01}v+\frac{1}{2}\left( c_{20}u^2+2c_{11}uv+c_{02}v^2 \right)+\cdots,
$$
$$
d=d_{0}+d_{10}u+d_{01}v+\frac{1}{2}\left( d_{20}u^2+2d_{11}uv+d_{02}v^2 \right)+\cdots.
$$
Since
$b_u=a_v$ and $d_u=c_v$, we have that  
$$
b_{10}=a_{01},\ d_{10}=c_{01},\ b_{20}=a_{11},\ b_{11}=a_{02}\ d_{20}=c_{11},\ d_{11}=c_{02},
$$
and so on.

\section{Ridges at Non-Discriminant Points}

In this section we obtain conditions in the jets of $a,b,c$ and $d$ for a non-discriminant point to be in the ridge set, for the
ridge set to be a regular curve, for the tangent to the ridge curve to coincide with the principal direction and for this 
contact to be of degree $2$. We start with a general lemma that will be useful also for points in the discriminant set. 

\subsection{Equations of the ridge set}

The eigenvalues of the shape operator are 
$$
\lambda_1=\frac{1}{2}\left( a+d+\sqrt{\delta} \right) \ \ {\rm and} \ \ \lambda_2=\frac{1}{2}\left( a+d-\sqrt{\delta} \right),
$$
where $\delta=(d-a)^2+4bc$.

\begin{lemma}\label{lemma:Useful}
Outside the discriminant curve, the equation of the ridge set $d\lambda_1(w_1)=0$ is $g_1=0$, $g_1=A\sqrt{\delta}+B$, where
$$
A=(a-d)a_u+cd_v+2ca_v+bc_u,
$$
$$
B=(a-d)^2a_u+(a-d)(2ca_v-cd_v+bc_u)+2c(2bd_u+ba_u+cb_v).
$$
\end{lemma}

\begin{proof}
We have that
$$
2d\lambda_1(w_1)=\frac{1}{2}(a-d+\sqrt{\delta})\left(a_u+d_u+\frac{\delta_u}{2\sqrt{\delta}}\right)+c\left(a_v+d_v+\frac{\delta_v}{2\sqrt{\delta}}   \right).
$$
Thus 
$$
8\sqrt{\delta}d\lambda_1(w_1)=(a-d+\sqrt{\delta})\left( 2\sqrt{\delta}(a_u+d_u)+\delta_u   \right)+2c\left(2\sqrt{\delta}(a_v+d_v)+\delta_v \right).
$$
The second member can be written as $A\sqrt{\delta}+B$, where
$$
A=4c(d_v+a_v)+2(a-d)(a_u+d_u)+2(a-d)(a_u-d_u)+4(b_uc+c_ub),
$$
$$
B=2c\left(2(a-d)(a_v-d_v)+4(b_vc+c_vb)\right)+2\delta(d_u+a_u)
$$
$$
+(a-d)\left(2(a-d)(a_u-d_u)+4(b_uc+c_ub)\right).
$$
This completes the proof.
\end{proof}

\subsection{Regularity of the ridge set}

In the complement of the discriminant set, the function $g_1$ defined in the above lemma is smooth. Similarly, the ridge curve 
$d\lambda_2(w_2)=0$ outside the discriminant curve can be written as $g_2=0$, for some smooth function $g_2$. 

Denote by $\mathcal{U}_2\subset\mathcal{U}_0$ the subset of line congruences such that the conditions $g_i=(g_i)_u=(g_i)_v=0$, $i=1,2$, do not occur simultaneously in the complement of the discriminant set. It is clear that, for line congruences in $\mathcal{U}_2$, the ridge set is a regular curve. 

Observe that the $3$ equations $g_i=(g_i)_u=(g_i)_v=0$  in $2$ variables $(u,v)$ involve up to third order derivatives of $\xi_1$ and $\xi_2$. Thus, by Thom's transversality theorem, the set $\mathcal{U}_2$ is open and dense in $\mathcal{U}_0$. From now on we shall assume that our line congruences belong to $\mathcal{U}_2\cap\mathcal{U}_1$.

\subsection{Principal direction tangent to the ridge curve }\label{sec:WTangentRidge}

At a point $p=(0,0)$ with two distinct principal directions,
we may assume that $(1,0)$ and $(0,1)$ are eigenvectors, which is equivalent to $b_0=0$ and $c_0=0$. 
Assume also that $p$ is not umbilic, which means that $a_0-d_0\neq 0$. We may assume, w.l.o.g., that $a_0-d_0>0$, which implies $\lambda_1(0,0)=a_0$, $w_1=(\lambda_1-d,c)$ is an eigenvector of $\lambda_1$ with $w_1(0,0)=(a_0-d_0,0)$
From Lemma \ref{lemma:Useful}, at the origin, $g_1=2(a_0-d_0)^2a_{10}$. Thus the origin is a $1$-ridge point if and only if $a_{10}=0$. In this section we shall assume that $a_{10}=0$.

\begin{lemma}\label{lemma:WTangentRidge}
The eigenvector $w_1$ is tangent to the ridge curve at the origin
if and only if $(a_0-d_0)a_{20} +3c_{10}a_{01}=0$.
\end{lemma}

\begin{proof}
Since, at the origin, $A=B=0$, we have that $(g_1)_u=A_u\sqrt{\delta}+B_u$.
Now straightforward calculations show that 
$$
(g_1)_u(0,0)=2(a_0-d_0)\left[ (a_0-d_0)a_{20} +3c_{10}a_{01} \right],
$$
which proves the lemma.
\end{proof}

\subsection{Contact of the ridge curve with the principal line}\label{sec:RidgeSecondOrder}

The principal line at the origin has a second order contact with the $u$-axis if and only if $c_{10}=0$. In fact, the angular coefficient of $w_1=(\lambda_1-d,c)$ satisfies
$$
\frac{c}{\lambda_1-d}=c_{10}u+O(u^2).
$$
By a change of variables of the form $V=v-\frac{c_{10}}{2}u^2$, we may assume that $c_{10}=0$. 

Consider a point where the principal direction is tangent to the ridge curve. By Lemma \ref{lemma:WTangentRidge}, since $c_{10}=0$, we have that $a_{20}=0$.

\begin{lemma}\label{lemma:WTangentRidge2}
At a non-discriminant ridge point, assume that $c_{10}=a_{20}=0$. Then the contact of the ridge curve with the principal line
is of order $2$ if and only if
\begin{equation}
4c_{20}a_{01} +(a_0-d_0)a_{30}\neq 0.
\end{equation}
\end{lemma}
\begin{proof}
Observe that, at the origin
$$
(g_1)_{uu}=A_{uu}(a_0-d_0)+B_{uu}=2a_{30}(a_0-d_0)^2+8(a_0-d_0)c_{20}a_{01},
$$
which proves the lemma.
\end{proof}

\section{Properties of the Discriminant Set}

Consider the local model and notation of Section \ref{sec:LocalModel}. We shall assume that the origin is a discriminant point with principal direction $(0,1)$. Then necessarily $a_0=d_0$, $b_0=0$ 
and $c_0\neq 0$. 

The eigenvector $(0,1)$ is tangent to the discriminant curve $\delta=0$ if and only if $b_{01}=0$. 
Points of the discriminant set with $b_{01}\neq 0$ will be called {\it regular} while points with $b_{01}=0$ will be called  {\it singular}.

At a singular point of the discriminant set, $\delta_v=0$ and $\delta_u=4c_0b_{10}$. Since $0$ is a regular value of $\delta$, we conclude that $b_{10}\neq 0$.

\subsection{Principal Directions}

In this section we describe the behaviour of the principal lines at the discriminant curve. Our approach follows \cite{Davydov}.
Let $p=\frac{du}{dv}$ and write Equation \eqref{eq:CurvatureCongruence2} as $F=0$, where
$$
F(u,v,p)=cp^2+(d-a)p-b.
$$
Then $F=0$ define a surface in the $(u,v,p)$ space. Since $\frac{\partial F}{\partial u}(0,0,0)=-b_{10}\neq 0$, 
we can write $F=0$ in a neighborhood of $(0,0,0)$ as $u=f(v,p)$, for some function $f$. Moreover, 
since $pdv=du=f_vdv+f_pdp$, we have the direction field given by
\begin{equation}\label{eq:EDBLift}
(f_v-p)dv+f_pdp=0.
\end{equation}
Note that at the origin, $f_p=\frac{d_0-a_0}{b_{10}}=0$ and $f_v=-\frac{b_{01}}{b_{10}}$. Thus $(p,v)=(0,0)$ is a singular point of Equation \eqref{eq:EDBLift} if and only if $b_{01}=0$, or equivalently, the discriminant point is singular. By \cite{Davydov}, the principal lines are cusps at regular points of the discriminant curve.

Assume now that the discriminant point is singular, i.e., $b_{01}=0$. The linear part at the origin of Equation \eqref{eq:EDBLift} is given by the matrix
\[
L=\frac{1}{a_{01}}
\left[
\begin{array}{cc}
d_{01}-a_{01} &  2c_{00}\\
b_{20} & 2a_{01}-d_{01}
\end{array}
\right]\ .
\]
Denote by the determinant of $L$
$$
\Delta_1=(d_{01}-a_{01})(2a_{01}-d_{01})-2c_0b_{20}.
$$
If $\Delta_1<0$, the origin is a folded saddle. If $\Delta_1>0$, the eigenvalues of $L$ are real (folded node) if $\Delta_2>0$ or complex (folded focus) if $\Delta_2<0$, where
$$
\Delta_2=(2d_{01}-3a_{01})^2+8c_0b_{02}.
$$

Denote by $\beta$, $p=\beta v$, the direction of the eigenvectors of $L$. The projections of the lines $p=\beta v$
are called the {\it separatrices} of the line congruence. 
Straighforward calculations show that, at the origin, 
\begin{equation}\label{eq:Separatriz}
2c_0\beta^2+(2d_{01}-3a_{01})\beta-b_{02}=0.
\end{equation}
The discriminant of this equation is exactly $\Delta_2$. Thus 
if $\Delta_2>0$, there are two separatrices while if $\Delta_2<0$ there are no separatrices.

\begin{example}\label{ex:Darbouxian}
Consider a Darbouxian umbilical point given by the BDE
$$
udu^2+mvdudv-udv^2=0,
$$
where $m$ is a parameter and the one-parameter unfolding 
$$
(u+s)du^2+mvdudv-(u-s)dv^2=0.
$$
This last BDE can be realized as the principal directions BDE of the line congruence with shape operator
$$
c=u+s,\ b=u-s,\ d=(m+1)v+1,\ a=v+1.
$$
The matrix $L$ at $(s,0,0)$ is then
\[
L=\left[
\begin{array}{cc}
m & 2s \\
0 & -m+1
\end{array}
\right]
\]
Thus the point is a folded saddle for $m<0$ and $m>1$, and a folded node for $0<m<1$. For $m>1$, this the only point of the discriminant set with tangent principal direction. For $m<1$, there are more $2$ such points, both folded saddles (see Figure \ref{Fig:Darboux}).





\begin{figure}[htb]
	\centering \hspace*{0.1 cm} \subfigure[Unfolded D1.]
{
	\includegraphics[
	trim=1cm 5cm 3cm 2.4cm, clip=true,width=.24
	\linewidth ]{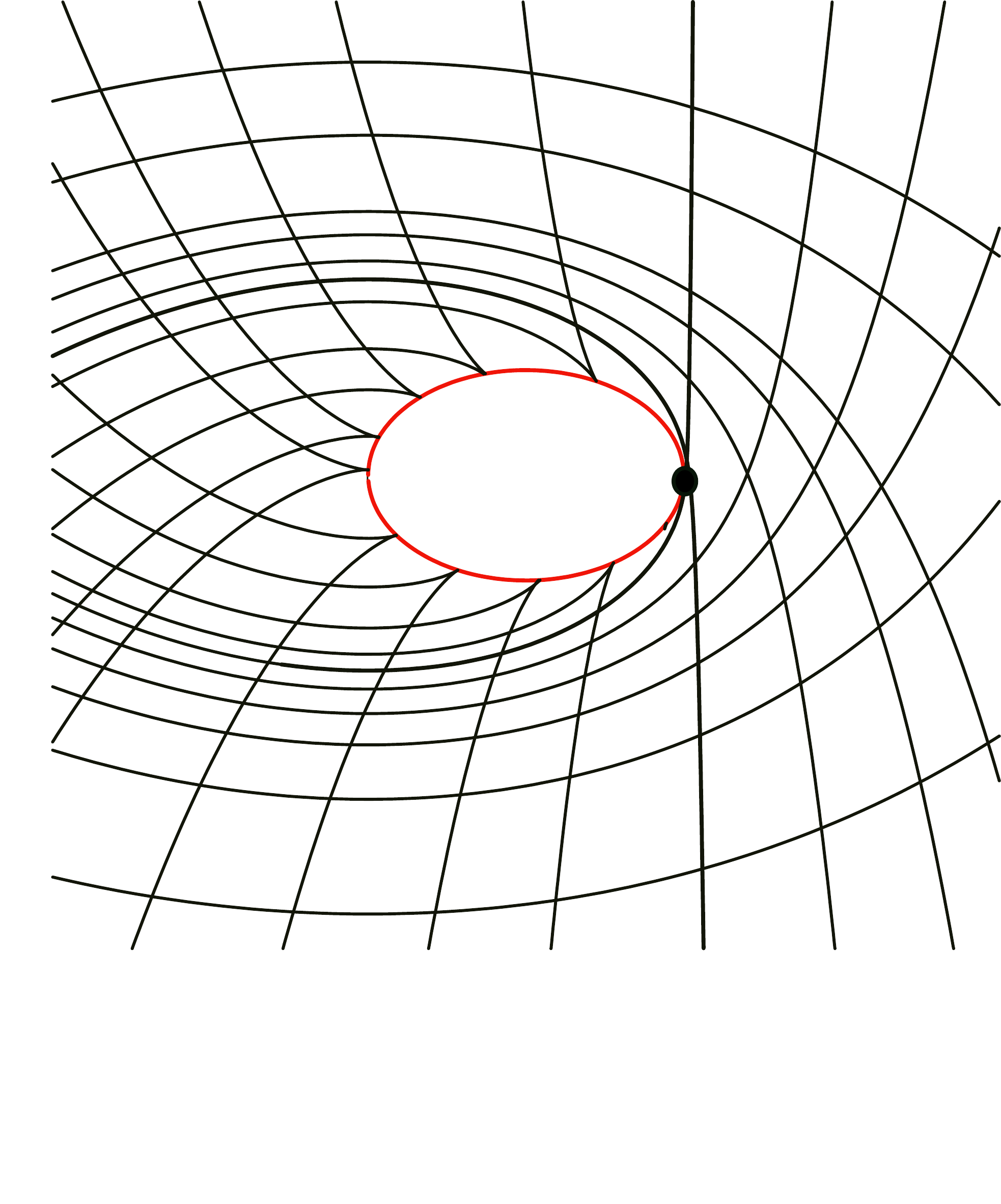}}\hspace*{0.5 cm}
\subfigure[ Unfolded D2.] 
{
\includegraphics[width=.25
\linewidth,clip =false]{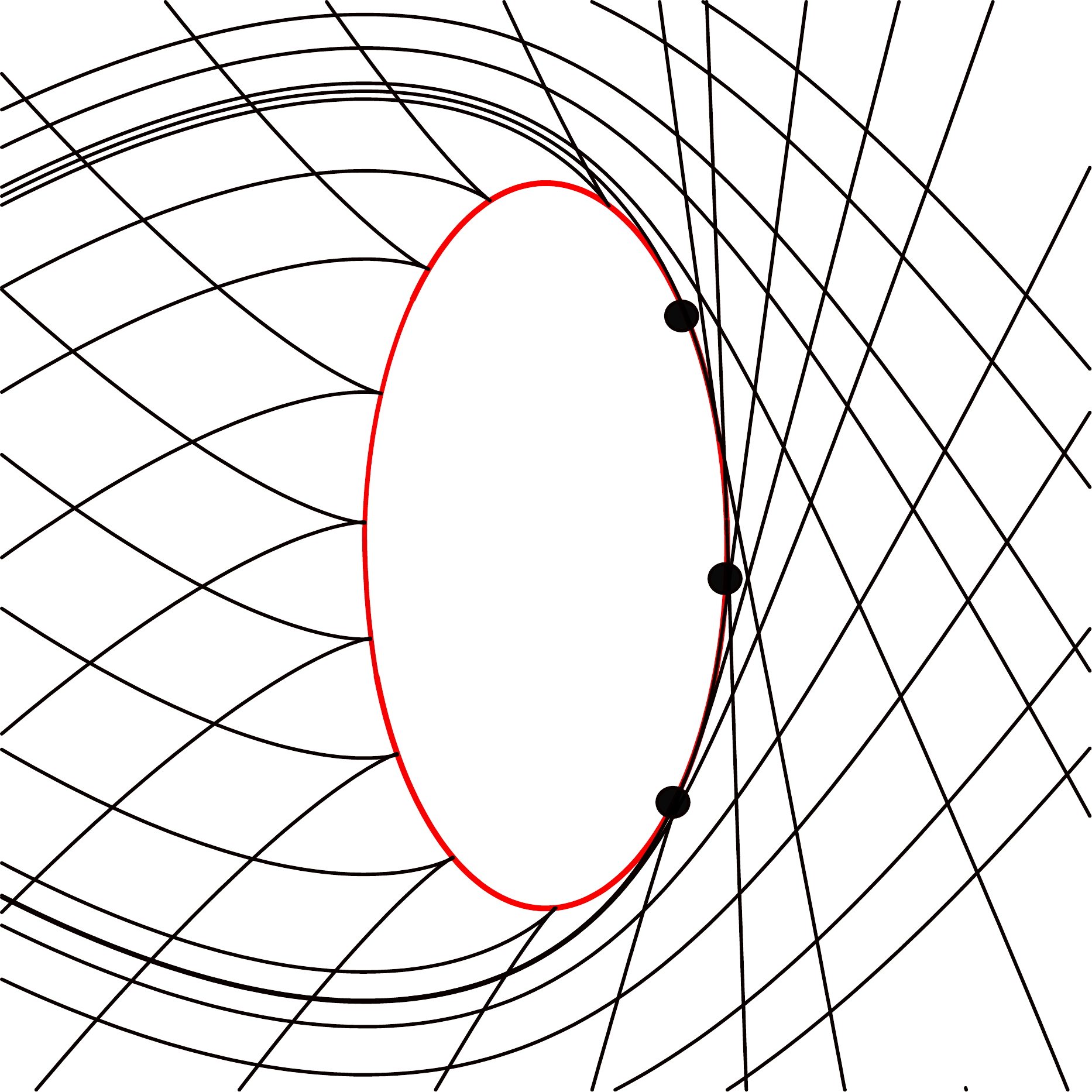}} \hspace*{0.5 cm}\subfigure
[Unfolded D3. ] 
{
\includegraphics[width=.24\linewidth,clip
=false]{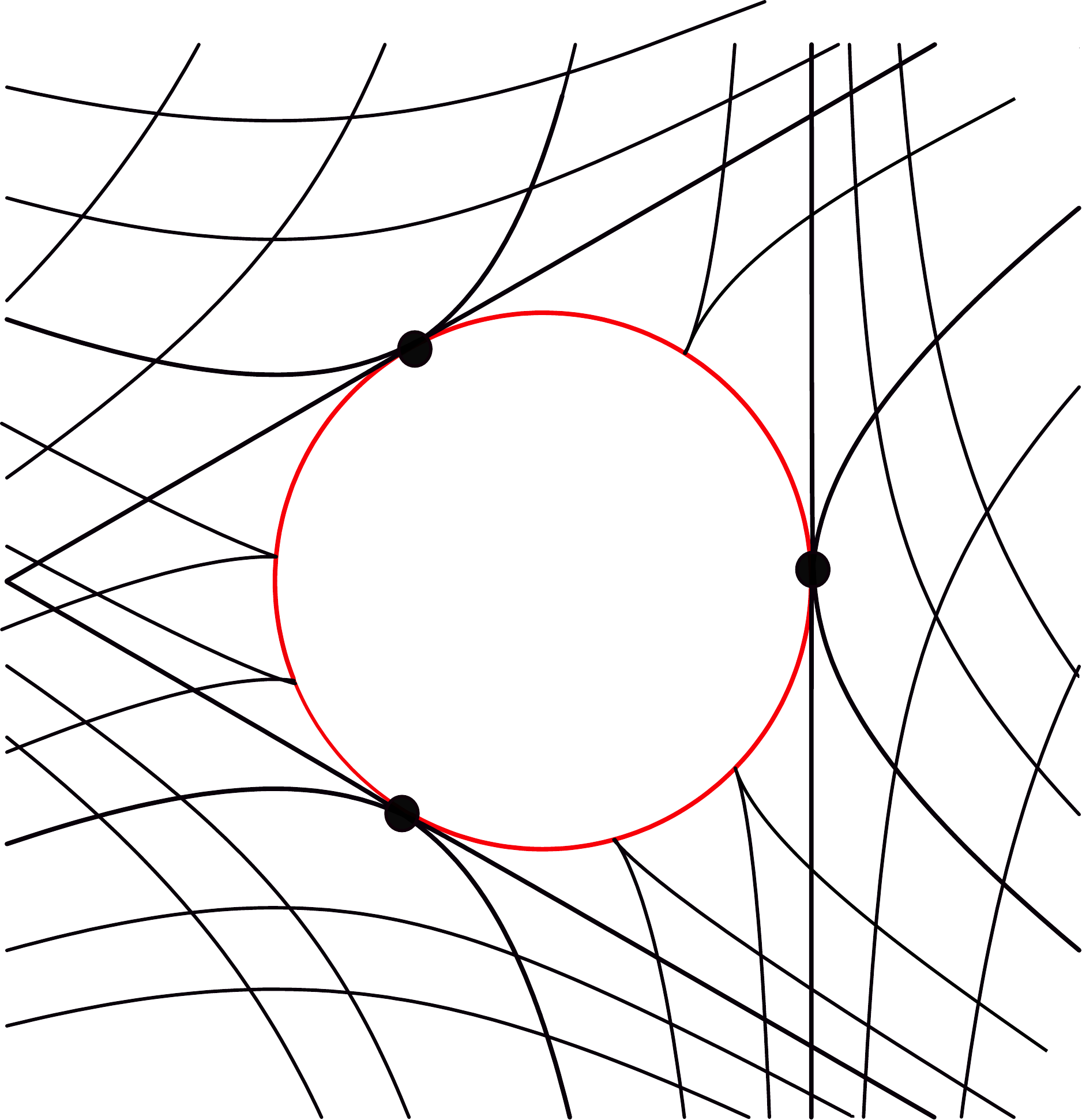}}\hspace*{0.1 cm}
\caption{Discriminant curve (red) and principal lines (black) of Example \ref{ex:Darbouxian} with $m<0$, $0<m<1$ and $m>1$, respectively.}
\label{Fig:Darboux}
\end{figure}

\end{example}

\subsection{Discriminant ridge points}

We consider the eigenvectors $w_i=(\lambda_i-d,c)$ associated to $\lambda_i$ with $w_i(0,0)=(0,c_0)$. We can then use 
Lemma \ref{lemma:Useful} to obtain that the origin is a ridge point if and only if $B=2c_0^2b_{01}=0$. Thus we can conclude that the ridge curve touches the discriminant curve only at singular points. In this section we shall assume $b_{01}=0$. 

Let $A$ and $B$ be given by Lemma \ref{lemma:Useful} and denote
$$
A_0=A(0,0)=c_0(2a_{01}+d_{01})
$$
and 
$$
B_1=B_v(0,0)=c_0\left(a_{01}(a_{01}-d_{01})+ (a_{01}-d_{01})^2+2c_0b_{02}\right).
$$
Denote also
$$
Q=\frac{1}{4a_{01}c_0}\delta_{vv}(0,0)=\frac{1}{2a_{01}c_0}\left( (a_{01}-d_{01})^2+2c_0b_{02} \right).
$$

\begin{Proposition}\label{prop:EqRidgeDiscriminant}
Assume that $A_0\neq 0$.
Then the ridge curve can be written locally as a graph $u=u(v)$ of the form 
$$
u(v)= \alpha\frac{v^2}{2}+O(v^3),
$$
where
$$
\alpha= \frac{B_1^2}{2c_0a_{01}A_0^2}-Q.
$$
Moreover, the part $v<0$ corresponds to one eigenvalue while the part $v>0$ corresponds to the other eigenvalue.  
\end{Proposition}

\begin{proof}
Assume w.l.o.g. that $A_0>0$ and $B_1>0$. Then the condition for ridge point associated with $\lambda_1$, $A\sqrt{\delta}+B=0$, can be written as $G=0$, where
$$
G=A^2\delta-B^2
$$
with $v<0$. A similar calculation shows that the condition for ridge point associated with $\lambda_2$ is given by $G=0$
with $v>0$. 

Observe that, at the origin 
$$
G_u=A_0^2\delta_u(0,0)=4A_0^2a_{01}c_0\neq 0
$$
and $G_v=0$. Moreover
$$
G_{vv}=A_0^2\delta_{vv}(0,0)-2B_1^2=4a_{01}c_0A_0^2Q-2B_1^2,
$$
thus proving the proposition.
\end{proof}

\subsection{Contact between a ridge curve and a separatrix}

\begin{Proposition}\label{prop:RidgeSeparatrix2}
Assume that $\Delta_1\neq 0$.
Then the contact between the ridge curve and a separatrix is of order $\geq 3$ if and only if
$$
3a_{01}d_{01}=c_0b_{02}.
$$
\end{Proposition}

\begin{proof}
If we substitue $\alpha$ given by Proposition \ref{prop:EqRidgeDiscriminant} in Equation \eqref{eq:Separatriz} we obtain 
$$
\left( 3a_{01}d_{01}-c_0b_{02}  \right)\left(2a_{01}^2+4a_{01}d_{01}-b_{02}c_0\right)\Delta_1^2=0.
$$
To conclude the proof we shall verify that the first factor must be zero.  
Observe that along the ridge we have 
$$
\delta=2a_{01}(c_0\alpha+Q)v^2=\frac{B_1^2}{A_0^2}v^2.
$$
Since $v<0$, the slope of the eigenvector $w_1$ is then
$$
\frac{1}{2c}\left( (a-d) -\frac{B_1}{A_0}v \right)=\frac{1}{2c_0}\left( (a_{01}-d_{01})-\frac{B_1}{A_0}\right)v+O(v^2).
$$
Hence the contact of the ridge curve with a principal line is of order greater than $2$ if and only if 
$$
(a_{01}-d_{01})-\frac{B_1}{A_0}-2c_0\alpha=0.
$$
This equation can be writen as 
$\left( 3a_{01}d_{01}-c_0b_{02}  \right)\Delta_1=0$,
thus proving the proposition.
\end{proof}

\begin{figure}[htb]
\centering
\includegraphics[width=0.3\linewidth]{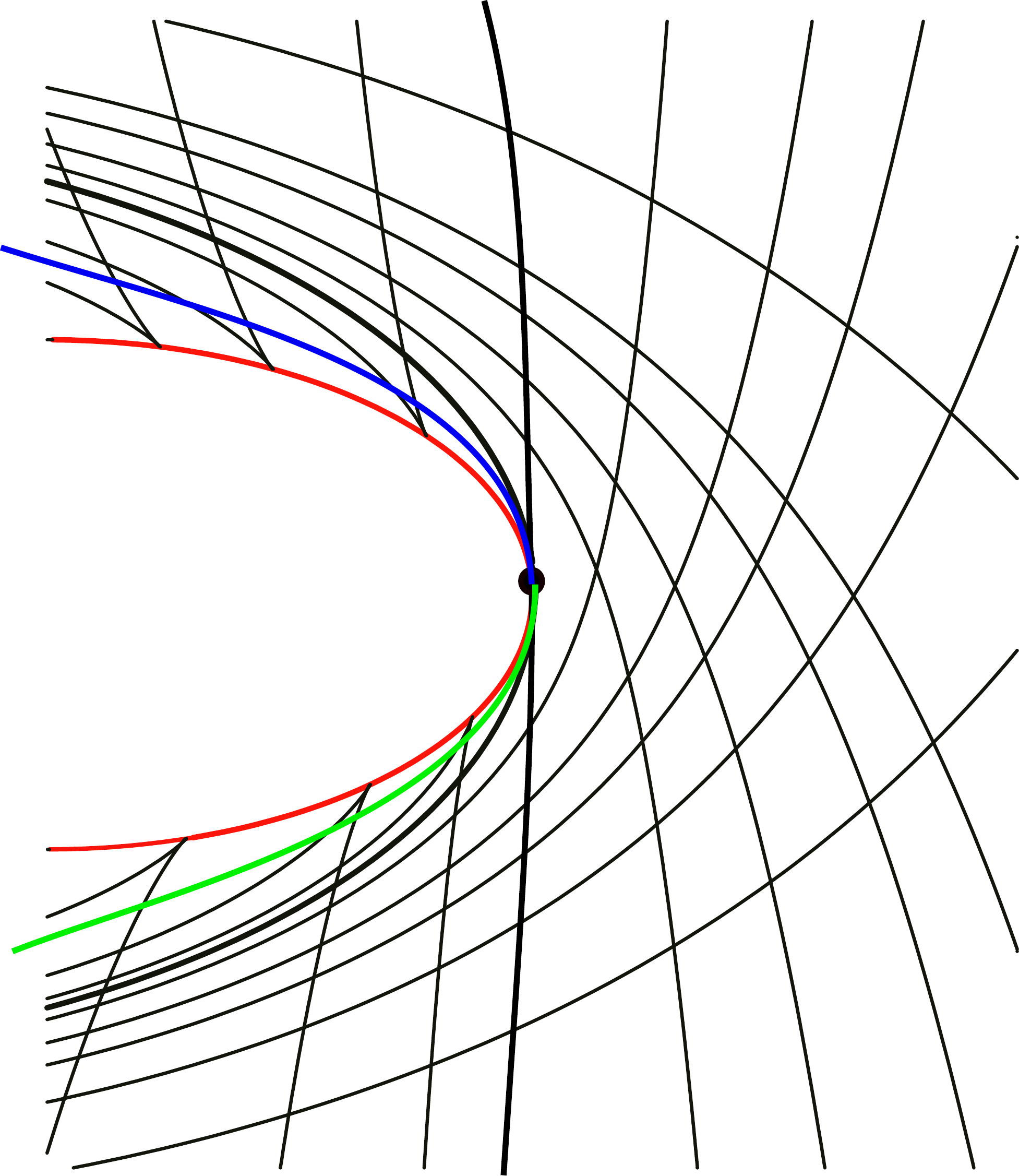}
 \caption{Discriminant curve (red), separatrices (thick black), and ridges (blue and green) at a singular discriminant point.}
\label{fig:DiscriminantPoint}
\end{figure}

In Figure \ref{fig:DiscriminantPoint}, one can see an illustration of the ridge curve making a contact of order $2$  
with the discriminant curve and the separatrices.

\subsection{Focal set at the discriminant curve}

In this section give a somewhat informal description of the focal set at the discriminant curve. The formal proofs are consequences of the singularity classification of the next sections.

At the discriminant curve, we write $\lambda_1=\lambda_2=\lambda=\frac{1}{2}(a+d)$. 
Equation \eqref{eq:Ft} with $\tau=0$ is written as
\begin{equation}\label{eq:Ft0}
D_XF_i=f_*\left( (I-\lambda^{-1}S)X    \right) -\lambda^{-2}D_X\lambda_i  \xi.
\end{equation}
Moreover, differentiating $2\lambda_i=a+d\pm\sqrt{\delta}$ in the direction $X\in\mathfrak{X}(U)$ we obtain
\begin{equation}\label{eq:dlambda}
2D_X\lambda_i=D_X(a+d) \pm \frac{D_X\delta}{2\sqrt{\delta}}.
\end{equation}

If $X$ is not tangent to the discriminant curve, then $D_X\delta\neq 0$ at a discriminant point. Thus 
Equation \eqref{eq:dlambda} implies that $D_X\lambda_i$ is is close to $\infty$ for a point close to the discriminant curve $\delta=0$. Then
Equation \eqref{eq:Ft0} implies that $D_XF_i$ is close to the $\xi$ direction for a point close to the discriminant curve.

If $X\in\mathfrak{X}(U)$ is tangent to the level sets of $\delta$ close to the discriminant curve, then
$D_X\lambda_i$ is close to $D_X(a+d)$. If moreover $X$ is not an eigenvector, $f_*\left( (I-\lambda^{-1}S)X    \right)$ is a non-zero vector in the image by $f_*$ of the eigenvector direction. Thus, by Equation \eqref{eq:Ft0}, $D_XF_i$ is a vector linearly independent from $\xi$. We conclude that, at a non-singular discriminant point, the focal set has a tangent plane generated by the image by $f_*$ of the eigenvector direction and $\xi$.

\begin{figure}[htb]
\centering
\includegraphics[width=0.21\linewidth]{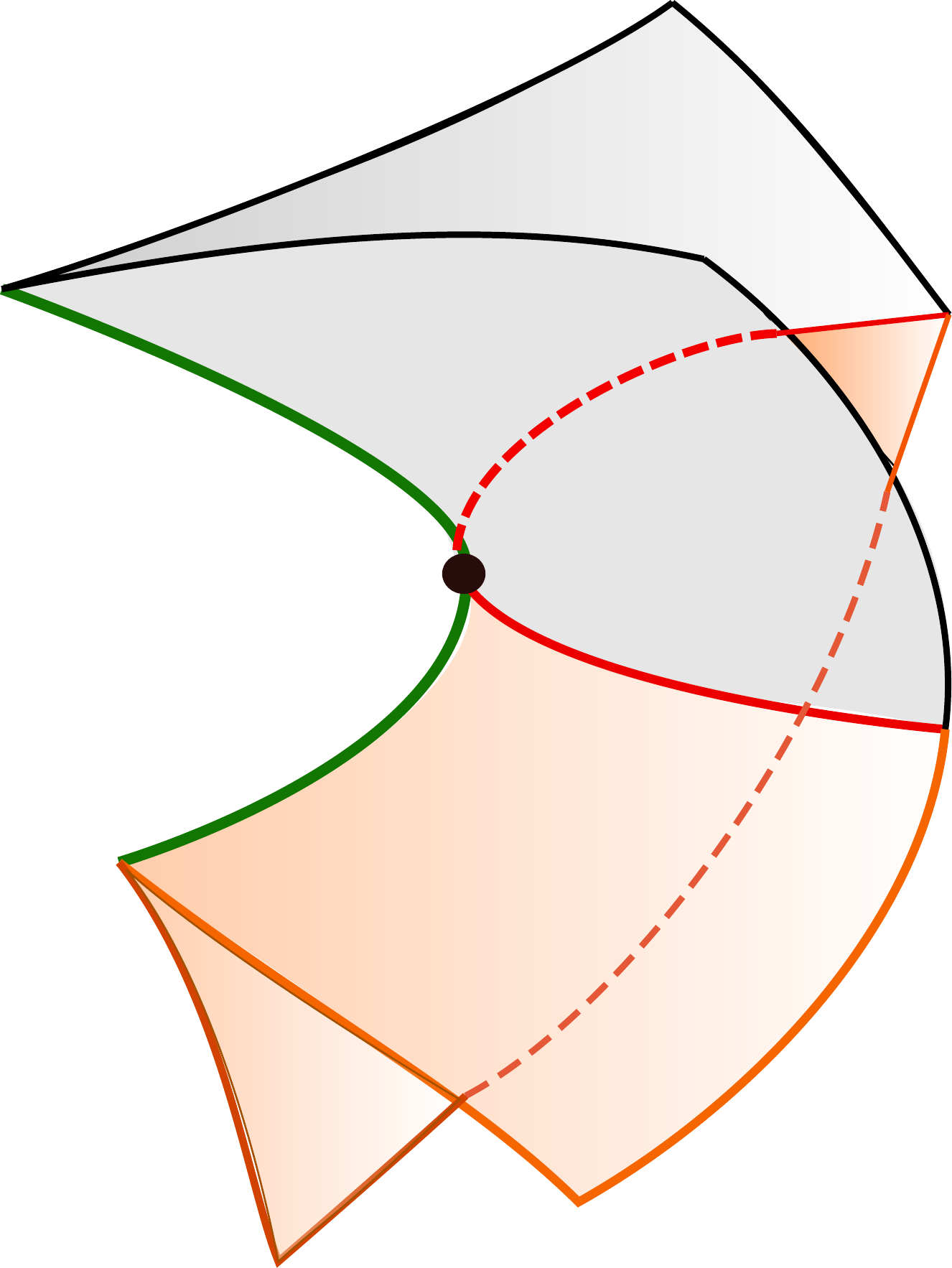}
 \caption{At a singular point, the image of the discriminant curve (red) separates the two branches of the focal set, while the image of the ridge curve (green) is the cuspidal edge.}
\label{fig:ModeloFocal}
\end{figure}

At a singular discriminant point, the tangent direction is also an eigenvector. Thus, by Equation \eqref{eq:Ft0}, the image of $DF_i$ is the $1$-dimensional subspace generated by $\xi$. 
In case $a_{01}+d_{01}\neq 0$, $D_XF_i$ is non-zero in the $\xi$-direction. Thus
the image of the discriminant curve by $F$ is regular with tangent vector in the $\xi$-direction. Moreover, from Proposition \ref{prop:EqRidgeDiscriminant}, the ridge changes color when touching the discriminant curve. Thus the image of the discriminant curve is passing from one focal surface $F_i$ to the other (see Figure \ref{fig:ModeloFocal}).
This configuration of focal sets and ridges appears also in the context of surfaces in Minkowski $3$-space (see \cite{Tari}, Figure 3).

\section{Equiaffine Pairs}

Given a pair $(f,\xi)$, we can write, for $X,Y\in\mathfrak{X}(U)$,
\begin{equation}\label{eq:DefineMetrica}
D_Xf_*Y=f_*(\nabla_XY)+h(X,Y)\xi,
\end{equation}
where $h$ is a metric. We say that the pair $(f,\xi)$ is {\it equiaffine}
if $h$ is non-degenerate and $\tau=0$, where $\tau$ is defined by Equation \eqref{eq:DefineShape2} (see \cite{Nomizu}). 

\subsection{Some examples}

\paragraph{Surfaces in $\R^3$} The normal vectors $\xi$ to a surfaces $f$ in the Euclidean $3$-space define a very well-known line congruence. Since the shape operator is self-adjoint, it admits two orthogonal eigenvectors at each point. The pair $(f,\xi)$ is equiaffine.

\paragraph{Surfaces in normed $3$-spaces}
Consider a norm in $\mathbb{R}^3$ whose unitary ball is a convex symmetric set $B$. Then the $B$-normal lines $\xi$ of a given surface $f$ defines an equiaffine pair. If we assume that $U$ is an ovaloid, then the corresponding shape operator admits two distinct eigenvalues at each point  (\cite{Balestro-Martini-Teixeira}).

\paragraph{Surfaces in Minkowski $3$-space}
Consider a surface in the Minkowski space $\mathbb{R}_1^3$. Outside the set of points with degenerate metric,
the unitary Minkowski normals define an equiaffine vector field. In the Riemannian part of the surface the discriminant curve is empty, but in the Lorentzian part it may be non-empty (see \cite{Tari}).

\paragraph{Blaschke affine vector field of a non-degenerate surface}
For a non-degenerate surface, the Blaschke vector field is an equiaffine vector field. For convex surfaces, the shape operator admits a pair of distinct eigenvectors at each point. For non-convex surfaces, the discriminant curve may be non-empty. The curvature lines for this type of line congruence was described in \cite{Barajas}.

\paragraph{Asymptotic lines of a surface in $\mathbb{R}^4$}
Every surface $G$ in $4$-space is locally the projective pedal of an equiaffine pair $(f,\xi)$, and the asymptotic lines
of $G$ correspond to the curvature lines of $(f,\xi)$ (\cite{Craizer-Garcia-2}).

\paragraph{Center symmetry sets}
Given an oval in $\mathbb{R}^3$, consider the two dimensional family of lines connecting points with parallel tangents. 
This constructiont is strongly related to equiaffine pairs.  In fact, the vector connecting the points with parallel tangents is an equiaffine vector field with respect to both surfaces. Conversely, if $(f,\xi)$ is an equiaffine pair, then the lines generated by $\xi$ is connecting points with parallel tangent planes of $f$ and $f+\xi$.
The focal set of this line congruence is called the {\it centre simmetry set} (\cite{Giblin-Zakalyukin}).

\subsection{Co-normal vector field and the support function}

Consider an eq\"uiaffine pair $(f,\xi)$. The co-normal vector field $\nu$ is defined by the conditions
\begin{equation*}
\nu\cdot f_u=\nu\cdot f_v=0,\ \ \nu\cdot\xi=1.
\end{equation*}
Since $(f,\xi)$ is eq\"uiaffine, we obtain $\nu\cdot\xi_u=\nu\cdot\xi_v=0$ and hence
$\nu_u\cdot \xi=\nu_v\cdot \xi=0.$
From $h_{11}=\nu\cdot f_{uu}$,  $h_{12}=\nu\cdot f_{uv}$ and $h_{22}=\nu\cdot f_{vv}$, we obtain that 
\begin{equation*}
\nu_u\cdot f_u=-h_{11}, \ \ \nu_v\cdot f_v=-h_{22},\ \ \nu_u\cdot f_v=\nu_v\cdot f_u=-h_{12}.
\end{equation*}
We can also write
\begin{equation}\label{eq:Shape}
-\left[
\begin{array}{cc}
\nu_{uu}\cdot \xi & \nu_{uv}\cdot \xi\\
\nu_{uv}\cdot \xi & \nu_{vv}\cdot \xi
\end{array}
\right]
= 
\left[
\begin{array}{cc}
\nu_u\cdot \xi_u & \nu_u\cdot \xi_v\\
\nu_v\cdot \xi_u & \nu_v\cdot \xi_v
\end{array}
\right]
= 
\left[
\begin{array}{cc}
h_{11} & h_{12}\\
h_{12} & h_{22}
\end{array}
\right]
\cdot
\left[
\begin{array}{cc}
a & b\\
c & d
\end{array}
\right] .
\end{equation}

The {\it support function} $\rho$ is defined as 
\begin{equation*}
\rho(u,v,Z)=\nu(u,v)\cdot (Z-f(u,v)),
\end{equation*}
where $\nu$ is the co-normal vector field and $Z=(x,y,z)$ is a point in $3$-space. Differentiating we obtain
\begin{equation*}
\rho_u=\nu_u\cdot (Z-f),\ \ \rho_v=\nu_v\cdot (X-f).
\end{equation*}
Thus $\rho_u=\rho_v=0$ if and only if $Z=f+t\xi$, for some $t=t(u,v)$. Differentiating once more we obtain
\begin{equation}\label{eq:D2Rho}
\rho_{uu}=\nu_{uu}\cdot(Z-f)+h_{11},\ \ \rho_{vv}=\nu_{vv}\cdot(Z-f)+h_{22},\ \ \rho_{uv}=\nu_{uv}\cdot(Z-f)+h_{12}.
\end{equation}
From Equations \eqref{eq:Shape}, we obtain that at $\rho_u=\rho_v=0$, 
\begin{equation}\label{eq:D2rho}
\left[
\begin{array}{cc}
\rho_{uu} & \rho_{uv}\\
\rho_{uv} & \rho_{vv}
\end{array}
\right]
=
\left[
\begin{array}{cc}
h_{11} & h_{12}\\
h_{12} & h_{22}
\end{array}
\right]
\cdot
\left(-t\left[
\begin{array}{cc}
a & b\\
c & d
\end{array}
\right] +I \right).
\end{equation}
The critical point $(u,v)$ of $\rho$ is thus degenerate if and only if
$t$ is equal to the inverse of an eigenvalue of the shape operator. 
We conclude that the focal set is the {\it bifurcation set} of the support function. 
The singularities of the line congruence are thus given by the degenerate singular points of $\rho$. Moreover, 
$A_i$-points of $\rho$, $i=2,3,4$, correspond respectively to fold, cusps and swallowtails of the line congruence.

\section{Line congruences as equiaffine pairs}\label{sec:EquiaffinePairs}

It is proved in (\cite{Giblin-Zakalyukin}) that, outside umbilical points, every line congruence can be obtained as the lines connecting points of two surfaces with parallel tangents. Inspired by the ideas of this paper, we prove in this section that, outside umbilical points, every line congruence can be represented by an equiaffine pair.

\subsection{Main result}

Let $f(u,v)=(u,v,0)$ and $\xi=(\xi_1,\xi_2,1)$. Then the shape operator is given by 
\begin{equation}\label{eq:ShapeOperator}
S=\left[
\begin{array}{cc}
a & b\\
c & d
\end{array}
\right] ,
\end{equation}
where $a=-(\xi_1)_u$, $b=-(\xi_1)_v$, $c=-(\xi_2)_u$ and $d=-(\xi_2)_v$.

\begin{Proposition}
In a neighborhood of a non-umbilical point of a line congruence, we can find an equiaffine pair representing the line congruence.
\end{Proposition}

\begin{proof}
The proof of this proposition will be divided in two steps: First we construct an equiaffine pair $(\tilde{f},\tilde{\xi})$, then prove that $h$ is non-degenerate separately for non-discriminant and discriminant points.

For the first set, let $p(u,v)$ be a solution of the partial differential equation
\begin{equation}\label{eq:EDP}
bp_{uu}+(d-a)p_{uv}-cp_{vv}=0.
\end{equation}
In general, this equation can be solved outside umbilical points (see \cite{Giblin-Zakalyukin}, last page).

Given $p$, we can find $q(u,v)$ satisfying 
$$
q_u=-(ap_u+cp_v),\ \ (\xi_3)_v=-(bp_u+dp_v).
$$
The co-normal $\nu$ is then defined as 
$$
\nu=\left( p_u, p_v, q-p_u\xi_1-p_v\xi_2 \right), 
$$
and the affine normal by $\tilde\xi=\frac{1}{q}\xi$. One can easily verify that 
$$
\nu\cdot\tilde\xi=1, \ \nu\cdot\tilde\xi_u=0,\ \nu\cdot\tilde\xi_v=0.
$$
Moreover define $k=-\frac{p}{q}$ and 
$$
\tilde{f}=\left( u+k\xi_1,v+k\xi_2,k\right).
$$
One can easily verify that $(\tilde{f},\tilde{\xi})$ is an equiaffine pair representing the line congruence $(f,\xi)$.

For the second step, consider a point outside the discriminant set. We may assume that, at the origin $d-a>0$ and $b=c=0$. In this case, we may choose a solution $p$ 
of Equation \eqref{eq:EDP} with $p=p_{u}=p_v=p_{uv}=0$ and $p_{uu}=p_{vv}=-1$ at the origin. Moreover we can choose
$q=1$ at the origin, which implies that $\nu=(0,0,1)$ and $k=0$ at the origin. Since $\xi=(0,0,1)$ we conclude that, at the origin $\tilde{f}=(0,0,0)$, $\tilde{f}_u=(1,0,0)$, $\tilde{f}_v=(0,1,0)$. Finally 
$\nu_u\cdot\tilde{f}_u=(\nu_1)_u=p_{uu}=-1$ at the origin. Similarly $\nu_u\cdot\tilde{f}_v=p_{uv}=0$ and $\nu_v\cdot\tilde{f}_v=p_{vv}=-1$. We conclude that $h_{11}=h_{22}=1$ and $h_{12}=0$, and so 
 the surface $\tilde{f}$ is convex and non-degenerate in a neighborhood of the origin.

At a point of the discriminant set, we may assume that, at the origin $a=d$ and $b=0$. We may choose a solution $p$ 
of Equation \eqref{eq:EDP} with $p=p_{u}=p_v=p_{uu}=p_{vv}=0$ and $p_{uv}=1$ at the origin. Moreover we can choose
$q=1$ at the origin, which implies that $\nu=(0,0,1)$ and $k=0$ at the origin. Since $\xi=(0,0,1)$ we conclude that, at the origin $\tilde{f}=(0,0,0)$, $\tilde{f}_u=(1,0,0)$, $\tilde{f}_v=(0,1,0)$. Finally 
$\nu_u\cdot\tilde{f}_u=(\nu_1)_u=p_{uu}=0$ at the origin. Similarly $\nu_u\cdot\tilde{f}_v=p_{uv}=1$ and $\nu_v\cdot\tilde{f}_v=p_{vv}=0$. We conclude that $h_{11}=h_{22}=0$ and $h_{12}=-1$, and so 
 the surface $\tilde{f}$ is not convex and non-degenerate in a neighborhood of the origin. 
\end{proof}

\paragraph{Parameters of the equiaffine pair}

The parameters $\tilde{a}, \tilde{b},\tilde{c}$ and $\tilde{d}$ are defined by the relations
$$
\tilde{\xi}_{u}=-\tilde{a}\tilde{f}_u-\tilde{c}\tilde{f}_v,\ \ \tilde{\xi}_{v}=-\tilde{b}\tilde{f}_u-\tilde{d}\tilde{f}_v.
$$

\begin{lemma}
The $1$-jet of $\tilde{a}$, $\tilde{b}$, $\tilde{c}$ and $\tilde{d}$ at the origin coincides with 
the $1$-jet of $a$, $b$, $c$ and $d$, respectively. 
\end{lemma}

\begin{proof}
At the origin, since $q_u=q_v=0$ and $q=1$, we obtain that $\tilde\xi_u=\xi_u$, $\tilde\xi_v=\xi_v$, $\tilde{f}_u=f_u$, $\tilde{f}_v=f_v$, which implies that 
$$
\tilde{a}_0=a_0,\ \tilde{b}_0=b_0\ \tilde{c}_0=c_0,\ \tilde{d}_0=d_0.
$$
Considering only the first two coordinates, at the origin $\tilde\xi_{uu}=\xi_{uu}$, $\tilde\xi_{uv}=\xi_{uv}$, 
$\tilde\xi_{vv}=\xi_{vv}$, $\tilde{f}_{uu}=f_{uu}$, $\tilde{f}_{uv}=f_{uv}$, $\tilde{f}_{vv}=f_{vv}$, which implies that
$$
\tilde{a}_{10}=a_{10}, \ \ \tilde{a}_{01}=a_{01},
$$
and the same hold for $b$, $c$ and $d$. 
\end{proof}

\subsection{Non-discriminant points}\label{sec:ParametersNonDiscriminant}

At non-discriminant points, we may assume that, at the origin $b=c=0$, $h_{11}=h_{22}=1$, $h_{12}=0$. 

\begin{lemma}
Assume $b_0=c_0=0$, $a_0-d_0\neq 0$. Then $\tilde{a}_{20}=a_{20}$ and $\tilde{c}_{20}=c_{20}$. Moreover if $a_{10}=c_{10}=0$, then $\tilde{a}_{30}=a_{30}$.
\end{lemma}

\begin{proof}
We have that
\begin{equation}\label{eq:Xi3u1}
\tilde{\xi}_{uuu}=-\tilde{a}_{uu}\tilde{f}_u-2\tilde{a}_u\tilde{f}_{uu}-\tilde{a}\tilde{f}_{uuu}-\tilde{c}_{uu}\tilde{f}_{v}-2\tilde{c}_u\tilde{f}_{uv}-\tilde{c}\tilde{f}_{uuv}.
\end{equation}
Considering only the first two coordinates at the origin, 
\begin{equation}\label{eq:TildeXi1}
\tilde\xi_{uuu}=-\tilde{a}_{uu}\tilde{f}_u-\tilde{a}\tilde{f}_{uuu}-\tilde{c}_{uu}\tilde{f}_v=\left(-\tilde{a}_{20}+3a_0^2,-\tilde{c}_{20}\right),
\end{equation}
since, at the origin,
$$
\tilde{f}_{uuu}=3k_{uu}(\xi_1)_u=-3a_0.
$$
On the other hand, 
\begin{equation}\label{eq:Xi3u2}
\tilde{\xi}_{uuu}=\frac{1}{q}\xi_{uuu}+3\left(\frac{1}{q}\right)_u\xi_{uu}+3\left( \frac{1}{q} \right)_{uu}\xi_u+\left( \frac{1}{q} \right)_{uuu}\xi.
\end{equation}
At the origin, considering the first two coordinates,
\begin{equation}\label{eq:TildeXi2}
\tilde{\xi}_{uuu}=\xi_{uuu}-3\frac{q_{uu}}{q^2}\xi_u=\left(-a_{20}+3a_0^2,-c_{20}\right).
\end{equation}
The first part of the lemma is proved by comparing Equations \eqref{eq:TildeXi1} and \eqref{eq:TildeXi2}. 

To prove the second part, differentiate Equations \eqref{eq:Xi3u1} and take the first coordinate to obtain
$$
(\tilde\xi_1)_{uuuu}=-\tilde{a}_{30}-a_0(\tilde{f_1})_{uuuu}.
$$
Since $(\tilde{f_1})_{uuuu}=4k_{uuu}(\xi_1)_u$ and we can choose $p_{uuu}(0,0)=0$, we obtain that $k_{uuu}=0$ and conclude that $(\tilde\xi_1)_{uuuu}=-\tilde{a}_{30}$. Now differentiate \eqref{eq:Xi3u2} with respect to $u$ at the origin and take the first coordinate to obtain
$$
(\tilde\xi_1)_{uuuu}=-a_{30}+4\left( \frac{1}{q}  \right)_{uuu}(\xi_1)_u
$$
Since $q_u=q_{uu}=q_{uuu}=0$ at the origin, we conclude that $(\tilde\xi_1)_{uuuu}=-a_{30}$ and the second part of the lemma is proved.
\end{proof}

\begin{lemma}
We have that, at the origin $\tilde{f}_{uv}=0$. Moreover,
$$
(h_{12})_u=\frac{a_{01}-c_{10}}{a_0-d_0}.
$$
\end{lemma}
\begin{proof}
Since $b=c=0$, we can take $p=p_u=p_v=p_{uv}=0$ and $p_{uu}=p_{vv}=-1$ at the origin. Differentiating the PDE 
with respect to $u$ at the origin we obtain
$$
-b_u+(d-a)p_{uuv}+c_u=0. 
$$
Since $c_v=d_u$, we conclude that, at the origin,
$$
p_{uuv}=-\frac{a_v-c_u}{a-d}
$$
Now one can check that $(\nu_u\cdot \tilde{f}_v)_u=p_{uuv}$ at the origin. Thus
$$
(h_{12})_u=-p_{uuv},
$$
proving the lemma.
\end{proof}

\subsection{Discriminant points}\label{sec:ParametersDiscriminant}

\begin{lemma}
Assume $a_0=d_0$, $b_0=0$, $c_0\neq 0$. Then $\tilde{b}_{02}=b_{02}$ and $\tilde{c}_{20}=c_{20}$.
\end{lemma}

\begin{proof}
Since
$$
\tilde{\xi}_{vv}=-\tilde{b}_v\tilde{f}_u-\tilde{b}\tilde{f}_{uv}-\tilde{d}_v\tilde{f}_{v}-\tilde{d}\tilde{f}_{vv},
$$
the first component at the origin of $\tilde{\xi}_{vvv}$ is 
\begin{equation}\label{eq:TildeXi3}
(\tilde{\xi_1})_{vvv}=-\tilde{b}_{02}-\tilde{d}\tilde{f}_{vvv}=-\tilde{b}_{02},
\end{equation}
since the first component of $\tilde{f}_{vvv}$ at the origin is zero. On the other hand, differentiating
$$
\tilde{\xi}_{vv}=\frac{1}{q}\xi_{vv}-2\frac{q_v}{q^2}\xi_v-\left( \frac{q_v}{q^2} \right)_v\xi
$$
with respect to $v$ we obtain, at the origin,
\begin{equation}\label{eq:TildeXi4}
(\tilde{\xi_1})_{vvv}=(\xi_1)_{vvv}=-b_{02}.
\end{equation}
Comparing Equations \eqref{eq:TildeXi3} and \eqref{eq:TildeXi4}, the lemma is proved.
\end{proof}

\begin{lemma}\label{lemma:DiscriminantMetric}
Assume $b_{01}=0$. At the origin, the following formulas hold:
$$
\nu_{vv}\cdot f_v=(\nu_2)_{vv}=k_{vvv}=p_{vvv}=\frac{d_{01}-a_{01}}{c_0}.
$$
$$
\nu_v\cdot f_{vv}=(\tilde{f}_1)_{vv}=0.
$$
\end{lemma}
\begin{proof}
From Equation \eqref{eq:EDP} we have that, at the origin, 
$$
(a_{01}-d_{01})h_{12}-c_0p_{vvv}=0
$$
Thus $p_{vvv}=\frac{d_{01}-a_{01}}{c_0}$. This proves the first claim. For the second claim, 
$$
(\tilde{f}_1)_{vv}=0
$$
since $k=k_v=k_{vv}=0$ at the origin. 
\end{proof}

\section{Singularities at Non-Discriminant Points}

Given $f=(u,v,0)$ and $\xi=(\xi_1,\xi_2,1)$, take the equiaffine pair $(\tilde{f},\tilde\xi)$ as in section $6$. The relevant parameters are not affected by the change from $(f,\xi)$ to $(\tilde{f},\tilde\xi)$, as proved in Section \ref{sec:ParametersNonDiscriminant}.

In this section we shall be dealing with the equiaffine pair $(\tilde{f},\tilde\xi)$, but in order to keep the notation shorter we
shall use the notation $(f,\xi)$ with the hope that this will cause no confusion.

\subsection{Singularities of type $A_2$ at non-discriminant points }

Differentiating
\begin{equation}\label{eq:Sup1}
\rho_{uu}=\nu_{uu}(Z-f)+h_{11},
\end{equation}
with respect to $u$ we obtain
\begin{equation}\label{eq:D3Rhou}
\rho_{uuu}=\nu_{uuu}(Z-f)-\nu_{uu}f_u+(h_{11})_u.
\end{equation}
Differentiating
\begin{equation}\label{eq:Sup2}
\nu_{uu}\cdot\xi=-h_{11}a-h_{12}c
\end{equation}
with respect to $u$ we obtain
\begin{equation}\label{eq:D3Etau}
\nu_{uuu}\cdot\xi-\nu_{uu}(af_u+cf_v)=-(h_{11})_ua-h_{11}a_u-(h_{12})_uc-h_{12}c_u.
\end{equation}
At the origin, we have
$$
\frac{1}{a_0}\nu_{uuu}\cdot\xi-\nu_{uu}f_u+(h_{11})_u=-\frac{a_{10}}{a_0}.
$$
We conclude that
$$
\rho_{uuu}=\frac{1}{a_0}\nu_{uuu}\cdot\xi-\nu_{uu}\cdot f_u+(h_{11})_u=-\frac{a_{10}}{a_0}.
$$
Thus if $a_{10}\neq 0$, the origin is a $A_2$-point of the support function. From Section \ref{sec:WTangentRidge}, we conclude the following:

\begin{corollary}
A point $p$ is a $A_2$ point of the support function if and only if $p$ is not a ridge point.
\end{corollary}

\subsection{Singularities of type $A_3$ at non-discriminant points}

From now on we shall assume $a_{10}=0$.

\begin{lemma}
At the origin,
$$
\rho_{uuv}=-\frac{c_{10}}{a_0}.
$$
\end{lemma}

\begin{proof}
Differentiating Equation \eqref{eq:D2Rho}
with respect to $u$ we obtain
$$
\rho_{uuv}=\nu_{uuv}(Z-f)-\nu_{uv}f_u+(h_{12})_{u}.
$$
At the origin, 
$$
\rho_{uuv}=\frac{1}{a_0}\nu_{uuv}\xi+\frac{1}{a_0}\nu_{uv}\xi_u+(h_{12})_{u}.
$$
Differentiating Equation \eqref{eq:D2Rho}
with respect to $u$ we obtain
$$
\nu_{uuv}\cdot\xi+\nu_{uv}\xi_u=-(h_{12})_ua-h_{12}a_u-(h_{22})_uc-h_{22}c_u.
$$
At the origin,
$$
\nu_{uuv}\cdot\xi+\nu_{uv}\xi_u+a(h_{12})_u=-c_u.
$$
Thus, at the origin, 
$$
\rho_{uvv}=-\frac{c_{10}}{a_0},
$$
thus proving the lemma.
\end{proof}

\begin{lemma}
At the origin,
$$
\rho_{uuuu}=-\frac{1}{a_0(a_0-d_0)}\left(  3(a_{01}-c_{10})c_{10}+(a_0-d_0)a_{20}  \right).
$$
\end{lemma}
\begin{proof}
Differentiating Equation \eqref{eq:D3Rhou}
with respect to $v$ we obtain
\begin{equation}\label{eq:Rho4u}
\rho_{uuuu}=\nu_{uuuu}(Z-f)-2\nu_{uuu}f_u-\nu_{uu}\cdot f_{uu}+(h_{11})_{uu}.
\end{equation}
At the origin
$$
\rho_{uuuu}=\frac{1}{a_0}\nu_{uuuu}\xi+\frac{2}{a_0}\nu_{uuu}\xi_u-\nu_{uu}\cdot f_{uu}+(h_{11})_{uu}.
$$
Since, at the origin
$$
\xi_{uu}=-c_uf_v-af_{uu},
$$
we obtain that 
$$
-\nu_{uu}\cdot f_{uu}=\frac{1}{a_0} \left(  \nu_{uu}\xi_{uu}+c_{10}\nu_{uu}\cdot f_v \right).
$$
On the other hand, since at the origin
$f_{uv}=0$, 
$$
\nu_{uu}\cdot f_v=-(h_{12})_u.
$$
We conclude that, at the origin
$$
\rho_{uuuu}=\frac{1}{a_0}\nu_{uuuu}\xi+\frac{2}{a_0}\nu_{uuu}\xi_u+\frac{1}{a_0}\nu_{uu}\cdot \xi_{uu}+(h_{22})_{uu}-\frac{c_{10}}{a_0}(h_{12})_u.
$$
Differentiating Equation \eqref{eq:D3Etau} with respect to $u$, at the origin, we obtain
$$
\nu_{uuuu}\cdot\xi+2\nu_{uuu}\cdot\xi_u+\nu_{uu}\cdot\xi_{uu}+a_0(h_{11})_{uu}=-2c_{10}(h_{12})_u-a_{20}.
$$
Thus
$$
\rho_{uuuu}=-\frac{1}{a_0}\left( 3c_{10}(h_{12})_u+a_{20}\right).
$$
Since $(h_{12})_u=\frac{a_{01}-c_{10}}{a_0-d_0}$, the lemma is proved.
\end{proof}

\begin{lemma}
The condition for an $A_3$ point of the support function is 
$$
3a_{01}c_{10}+(a_0-d_0)a_{20}\neq 0.
$$
\end{lemma}

\begin{proof}
The condition for $A_3$ point is $3\rho_{uuv}^2\neq \rho_{vv}\rho_{uuuu}$. 
\end{proof}

From Lemma \ref{lemma:WTangentRidge}, we conclude the following:

\begin{corollary}
A point $p$ is a $A_3$ point of the support function if and only if $p$ is a ridge point but the principal direction is not tangent to the ridge curve.
\end{corollary}

\subsection{Singularities of type $A_4$ at non-discriminant points}

In this section we shall assume that the principal direction is tangent to the ridge curve. In terms of jets, we have that $a_{10}=0$ and the condition of Lemma \eqref{lemma:WTangentRidge} holds. As mentioned in Section \ref{sec:RidgeSecondOrder}, we may in fact assume that
$c_{10}=a_{20}=0$. 

\begin{lemma}
Assume that $c_{10}=a_{20}=0$. Then the point is an $A_4$ point for $\rho$ if and only if
\begin{equation}
4a_{01}c_{20}+(a_0-d_0)a_{30}\neq 0.
\end{equation}
\end{lemma}

\begin{proof}
From Equation \eqref{eq:Rho4u} we obtain that, at the origin,
$$
\rho_{uuuuu}=\frac{1}{a_0}\nu_{uuuuu}\cdot\xi+\frac{3}{a_0}\nu_{uuuu}\cdot\xi_u-3\nu_{uuu}\cdot f_{uu}-\nu_{uu}\cdot f_{uuu}+(h_{11})_{uuu}.
$$
But 
$$
\nu_{uuu}\cdot f_{uu}=-\frac{1}{a_0}\nu_{uuu}\cdot\xi_{uu}, \ \ \nu_{uu}\cdot f_{uuu}=-\frac{1}{a_0}\nu_{uuu}\cdot\xi_{uuu}+\frac{c_{20}}{a_0}\nu_{uu}\cdot f_v.
$$
We conclude that
$$
\rho_{uuuuu}=\frac{1}{a_0}(\nu_{uu}\cdot\xi)_{uuu}+(h_{11})_{uuu}-\frac{c_{20}}{a_0}(h_{12})_u
$$
But
$$
(\nu_{uu}\cdot\xi)_{uuu}=-a_0(h_{11})_{uuu}-a_{30}-3(h_{12})_uc_{20}.
$$
Thus
$$
\rho_{uuuuu}=-\frac{a_{30}}{a_0}-4\frac{c_{20}}{a_0}\frac{a_{01}}{a_0-d_0}=-\frac{1}{a_0(a_0-d_0)}\left( 4a_{01}c_{20}+(a_0-d_0)a_{30} \right),
$$
thus proving the lemma.
\end{proof}

From Lemma \ref{lemma:WTangentRidge2}, we conclude the following:

\begin{corollary}
A point $p$ is a $A_4$ point of the support function if and only if $p$ is a ridge point, the principal direction is tangent to the ridge curve and the contact between the principal line and the ridge curve is exactly $2$. 
\end{corollary}

\section{Singularities at Discriminant Points}

Given $f=(u,v,0)$ and $\xi=(\xi_1,\xi_2,1)$, take the equiaffine pair $(\tilde{f},\tilde\xi)$ as in section $6$. The relevant parameters are not affected by the change from $(f,\xi)$ to $(\tilde{f},\tilde\xi)$, as proved in Section \ref{sec:ParametersDiscriminant}.

In this section we shall be dealing with the equiaffine pair $(\tilde{f},\tilde\xi)$, but in order to keep the notation shorter we
shall use the notation $(f,\xi)$ with the hope that this will cause no confusion.

\subsection{Singularities of type $A_2$ at discriminant points }

A singular point is a fold point of $G$ or a regular point of the focal set if and only $\rho_{vvv}\neq 0$ at this point. We shall prove now that this is equivalent to $b_{01}\neq 0$, which means that the eigenvector is not tangent to the double eigenvalues curve.

Differentiating
\begin{equation}\label{eq:Sup1}
\rho_{vv}=\nu_{vv}(Z-f)+h_{22},
\end{equation}
with respect to $v$ we obtain
\begin{equation}\label{eq:Sup12}
\rho_{vvv}=\nu_{vvv}(Z-f)-\nu_{vv}f_v+(h_{22})_v.
\end{equation}
At the double eigenvalues curve, we have $a=d$, $b=0$. We may also assume that, at this point,
$h_{11}=h_{22}=0$, $h_{12}=1$. Differentiating
\begin{equation}\label{eq:Sup2}
\nu_{vv}\cdot\xi=-h_{12}b-h_{22}d
\end{equation}
with respect to $v$ we obtain
\begin{equation}\label{eq:Sup22}
\nu_{vvv}\cdot\xi-\nu_{vv}(bf_u+df_v)=-(h_{12})_vb-h_{12}b_v-(h_{22})_vd-h_{22}d_v.
\end{equation}
At the point, we have
$$
\rho_{vvv}=\frac{1}{d_0}\nu_{vvv}\cdot\xi-\nu_{vv}\cdot f_v+(h_{22})_v=-\frac{1}{d_0} b_{01}.
$$

We have thus proved the following:

\begin{Proposition}
A discriminant point $p$ is an $A_2$ point if and only if it is regular. 
\end{Proposition}

\subsection{Singularities of type $A_3$ at discriminant points}

We now describe conditions for a singular point of the discriminant to correspond to a cuspidal edge of the focal set.
We shall assume $b_v=0$.

\begin{lemma}
At the origin,
$$
\rho_{uvv}=-\frac{b_{10}+d_{01}-a_{01}}{d_0}.
$$
\end{lemma}

\begin{proof}
Differentiating Equation \eqref{eq:Sup1}
with respect to $u$ we obtain
$$
\rho_{uvv}=\nu_{uvv}(Z-f)-\nu_{vv}f_u+(h_{22})_{u}.
$$
Differentiating Equation \eqref{eq:Sup2}
with respect to $u$ we obtain
$$
\nu_{uvv}\cdot\xi-\nu_{vv}(af_u+cf_v)=-(h_{12})_ub-h_{12}b_u-(h_{22})_ud-h_{22}d_u.
$$
At the origin,
$$
\rho_{uvv}=\frac{1}{d_0}\nu_{uvv}\cdot\xi-\nu_{vv}\cdot f_u+(h_{22})_u=-\frac{b_u}{d_0}+\frac{c_0}{d_0}\nu_{vv}\cdot f_v.
$$
From Lemma \ref{lemma:DiscriminantMetric}, 
$$
\rho_{uvv}=-\frac{b_{10}+d_{01}-a_{01}}{d_0},
$$
thus proving the lemma.
\end{proof}

\begin{lemma}
At the origin,
$$
\rho_{vvvv}=\frac{1}{c_0d_0}\left(  3(a_{01}-d_{01})d_{01}-c_0b_{02}  \right).
$$
\end{lemma}
\begin{proof}
Differentiating Equation \eqref{eq:Sup12}
with respect to $v$ we obtain
$$
\rho_{vvvv}=\nu_{vvvv}(Z-f)-2\nu_{vvv}f_v-\nu_{vv}\cdot f_{vv}+(h_{22})_{vv}.
$$
Differentiating Equation \eqref{eq:Sup22} with respect to $v$ leads to, at the origin
$$
\nu_{vvvv}\cdot\xi+2\nu_{vvv}\cdot\xi_v+\nu_{vv}\cdot\xi_{vv}+d_0(h_{22})_{vv}=-b_{vv}+2(h_{22})_vd_v.
$$
Since, at the origin,
$$
\xi_v=-d_0f_v,\ \ \xi_{vv}=-d_vf_{v}-d_0f_{vv},
$$
we obtain
$$
\nu_{vv}\cdot\xi_{vv}=-d_0\nu_{vv}\cdot f_{vv}-d_v\nu_{vv}\cdot f_v
$$
and so
$$
d_0\rho_{vvvv}=d_v\nu_{vv}\cdot f_v-b_{vv}-2(h_{22})_vd_v.
$$
But differentiating $h_{22}=\nu\cdot f_{vv}$ we obtain, at the origin, since $f_{vv}=0$, 
$$
(h_{22})_v=\nu\cdot f_{vvv}=g_{vvv}.
$$
Thus
$$
d_v\nu_{vv}\cdot f_v-b_{vv}-2(h_{22})_vd_v=-3g_{vvv}d_v-b_{vv}.
$$
We conclude that
$$
\rho_{vvvv}=\frac{1}{d_0}\left(  \frac{3(a_{01}-d_{01})d_{01}}{c_0}-b_{02}  \right),
$$
thus proving the lemma.
\end{proof}

\begin{Proposition}
A singular point of the discriminant set is an $A_3$-point of the support function if and only if
$$
3d_{01}a_{01}\neq c_0b_{02}.
$$
\end{Proposition}

\begin{proof}
The condition for cuspidal edge is that
$$
\frac{1}{2}\rho_{uu}u^2+\frac{1}{2}\rho_{uvv}uv^2+\frac{1}{24}\rho_{vvvv}v^4
$$
is not a perfect square, which is equivalent to
$$
3\rho_{uvv}^2\neq \rho_{uu}\rho_{vvvv}.
$$
Straightforward calculations complete the proof.
\end{proof}

From Proposition \ref{prop:RidgeSeparatrix2} we obtain the our main result:

\begin{corollary}
A singular point of the discriminant curve is a cuspidal edge if and only if the contact between the ridges and the separatrix
is exactly $2$.
\end{corollary}

\end{document}